\documentclass[final,leqno,letterpaper]{article}
\usepackage{amsmath}
\usepackage{hyperref}
\usepackage{graphicx}

\hypersetup{
    pdftitle={Solving the $3\times3$ Real Symmetric Eigenproblem},
    pdfauthor={Carlos F. Borges},
    pdfkeywords={arrow matrix, symmetric eignproblem, zero finder}
    }

\title{Solving the $3\times3$ Real Symmetric Eigenproblem}

\author{Carlos F. Borges\footnotemark[2]}

\def\bne{\begin{equation}}
\def\ene{\end{equation}}
\def\bmat{\left[ \begin{array}}
\def\endbmat{\end{array} \right]}

\def\b0{\mbox{\bf 0}}

\def\be{\mbox{\bf e}}

\def\bu{\mbox{\bf u}}

\def\fp{f^{\prime}}

\def\l{\lambda}

\def\Ab{\bar{A}}
\def\ab{\bar{\alpha}}
\def\gb{\bar{\gamma}}

\def\matlab{\textsc{Matlab }}

\begin{document}

\maketitle

\renewcommand{\thefootnote}{\fnsymbol{footnote}}

\footnotetext[2]{Authors address: Department of Applied Mathematics, Naval Postgraduate School, Monterey, CA 93943.  Email: borges@nps.edu}

\begin{abstract}
We develop an algorithm solving the $3\times3$ real symmetric eigenproblem. This is a common problem and in certain applications it must be solved many thousands of times, see for example \cite{tripref} where each element in a finite element grid generates one. Because of this it is useful to have a tailored method that is easily coded and compact. Furthermore, the method described is fully compatible with development as a GPU based code that would allow the simultaneous solution of a large number of these small eigenproblems.
\end{abstract}

\section{Reduction to arrow form}

The traditional first step in solving any real symmetric eigenproblem is to use unitary similarity transformations to reduce the matrix to tridiagonal (Hessenberg) form. For a $3\times3$ this would typically involve a Givens rotation that eliminates the $(1,3)$ and $(3,1)$ elements. We will stray from this approach and instead rely on a unitary similarity that eliminates the $(1,2)$ and $(2,1)$ elements and transforms our original matrix into an {\em ordered} $3\times3$ real symmetric arrow matrix.\footnote{I use an arrow structure rather than a tridiagonal for pedagogical reasons as I believe the derivations are easier to follow in this form. The two are equivalent under a simple permutation similarity.} Specifically, a matrix of the form
$$
A = 
\bmat{ccc}
\alpha_1 & 0 & \beta_1 \\
0 & \alpha_2 & \beta_2 \\
\beta_1 & \beta_2 & \gamma \\
\endbmat
$$
with $\alpha_1 \geq \alpha_2$. 

Although such a similarity can be constructed in several ways, we will use a Jacobi rotation (in effect, the eigenvectors of the principal $2\times2$) operating in the $(1,2)$-plane. The rotation can always be constructed so that the resulting arrow will have the elements of the shaft properly ordered. It is important to construct the Jacobi rotation carefully and we do so using the algorithm presented in \cite{borgesjacobi} which is demonstrably superior to the standard approach. 

\subsection{Deflating the Arrow}

We note that if any $\beta_j = 0$ then it is possible to set $\l_j = \alpha_j$ and deflate the matrix since $\be_j$ is clearly an eigenvector \cite{wilkinson}.  This is known as $\beta$-deflation. A second type of deflation occurs if $\alpha_1 = \alpha_2$, in that case we can apply a $2\times 2$ rotation similarity transformation in the $(1,2)$-plane that takes $\beta_1$ to zero and creates a $\beta$-deflation. This is known as a combo-deflation.

Exact deflations are theoretically easy to handle, however, before proceeding we must address the issue of {\em numerical deflation} which happens whenever {\em tiny} changes in the matrix can lead to deflation. Because the matrix is only $3\times3$ we can deal with the problem in a very direct manner. We begin by constructing a Givens rotation in the $(1,2)$-plane that takes $\beta_1$ to zero. The first step is to compute $h = \sqrt{\beta_1^2+\beta_2^2}$. if $h = 0$ then the matrix is diagonal and we are finished. If not, then we let
$$
G = \bmat{c c c} \beta_2/h & -\beta_1/h & 0 \\ \beta_1/h & \beta_2/h & 0 \\ 0 & 0 & 1 \endbmat.
$$
Then
$$
G A G^T = 
\bmat{ccc}
\frac{\alpha_1 \beta_2^2+\alpha_2 \beta_1^2}{h^2} & \alpha\frac{\beta_1\beta_2}{h^2} & 0 \\
\alpha\frac{\beta_1\beta_2}{h^2} & \frac{\alpha_1 \beta_1^2+\alpha_2 \beta_2^2}{h^2} & h \\
0 & h & \gamma \\
\endbmat
$$
where $\alpha = \alpha_1-\alpha_2$.

As the matrix is now tridiagonal we can adapt the deflation condition that is used in EISPACK (see \cite{gvl} pp.352-353) to our cause which would have us deflate if
$$
|\alpha\beta_1\beta_2| \leq C \epsilon \left( |\alpha_1 \beta_2^2+\alpha_2 \beta_1^2| + |\alpha_1 \beta_1^2+\alpha_2 \beta_2^2 | \right)
$$
for some constant $C$. We can simplify this condition by noting that
$$
|\alpha_1+\alpha_2| h^2\leq |\alpha_1 \beta_2^2+\alpha_2 \beta_1^2| + |\alpha_1 \beta_1^2+\alpha_2 \beta_2^2 | 
$$
by the triangle inequality. This leads to a simpler but slightly more restrictive test where if
$$
|\alpha\beta_1\beta_2| \leq C \epsilon  |\alpha_1+\alpha_2| h^2 
$$
we invoke the Wielandt-Hoffman theorem and ignore it yielding a deflation where we accept
$$
\frac{\alpha_1 \beta_2^2+\alpha_2 \beta_1^2}{h^2}
$$
as an eigenvalue. The rest of the spectrum can be recovered by solving the $2\times2$ eigenproblem for
$$
\bmat{cc}
\frac{\alpha_1 \beta_1^2+\alpha_2 \beta_2^2}{h^2} & h \\
 h & \gamma \\
\endbmat
$$

If the matrix does not deflate then we can assume that it is both ordered and {\em reduced}, that is the strict inequality $\alpha_1 > \alpha_2$ holds, and further $\beta_j \neq 0$ for $i=1,2$.

\subsection{Solving the Eigenproblem for a $3\times3$ Ordered and Reduced Arrow}

The interlacing property for real symmetric matrices combined with the fact that $\alpha_1 > \alpha_2$ implies that there is a rightmost eigenvalue $\l_1 > \alpha_1$ of multiplicity one with associated eigenvector $\bu_1$. If we shift $A$ by $\alpha_1$ we find that  $\mu = \l_1-\alpha_1$ is the only positive eigenvalue of
$$
A-\alpha_1 I = \bmat{ccc}
0 & 0 & \beta_1 \\
0 & -\ab & \beta_2 \\
\beta_1 & \beta_2 & \gamma-\alpha_1 \\
\endbmat
$$
where $\ab = \alpha_1-\alpha_2$. We shall call a symmetric matrix of this form (all zeros above the main counterdiagonal and a strictly negative middle element) a {fully reduced arrow}. The eigenvector associated with $\mu$ is identical with $\bu_1$ and it is easily verified that it is given by
$$
\bu_1 = \bmat{c} \beta_1 (\mu+\ab) \\  \beta_2 \mu \\ \mu(\mu+\ab) \endbmat
$$

Observe next that
$$
P(\alpha_2 I - A)P = \bmat{ccc}
0 & 0 & -\beta_2 \\
0 & -\ab & -\beta_1 \\
-\beta_2 & -\beta_1 & \alpha_2-\gamma \\
\endbmat
$$
where $P$ is the simple permutation  that swaps the first and second rows, is also a fully reduced arrow and that its only positive eigenvalue, $\nu$, satisfies $\nu  = \alpha_2 - \l_3$. One can verify that its associated eigenvector is
$$
\bmat{c}  \beta_2 (\nu+\ab) \\  \beta_1 \nu \\ -\nu(\nu+\ab) \endbmat
$$
and hence we can recover that corresponding leftmost eigenvector of $A$ by applying the permutation $P$ to get
$$
\bu_3 = \bmat{c} \beta_1 \nu  \\  \beta_2 (\nu+\ab) \\ -\nu(\nu+\ab) \endbmat
$$

Because the eigenvectors must be orthogonal we can compute the middle eigenvector directly by taking the cross product $\bu_1\times\bu_3$. This is better accomplished after a bit of algebra by using the form:
$$
\bu_2 = \bmat{c} 
-\beta_2 \mu (\nu+\ab)\\ \quad\beta_1\nu(\mu+\ab)\\ \beta_1\beta_2\ab
\endbmat
$$
It is worth noting that, discounting errors accrued in finding $\mu$ and $\nu$, there is no cancellation in the construction of any of the eigenvectors beyond the benign cancellation in computing $\ab$.

Finally, we can use the trace identity 
$$\mbox{tr}(A) = \l_1 + \l_2 + \l_3$$
to compute the corresponding eigenvalue. A bit of algebra yields
$$
\l_2 = \nu-\mu+\gamma
$$

In section \ref{algsec} we develop an algorithm that computes the dominant eigenvalue of a fully reduced arrow.

\subsection{Orthogonality of the Computed Eigenvectors}

In this section we show that the computed eigenvectors will be numerically orthogonal provided that
\begin{eqnarray}
\hat{\mu} & = & \mu (1+\delta_1) \nonumber \\
\hat{\nu} & = & \nu (1+\delta_2) 
\label{munu}
\end{eqnarray}
where $|\delta_i| < \epsilon$ for some small positive number $\epsilon$.

We begin by noting that since $\bu_2$ is a cross product it will be necessarily be numerically orthogonal to $\bu_1$ and $\bu_3$ provided that these two are themselves numerically orthogonal to each other (see \cite{kahan}) and so we only need demonstrate that. To that end we note that
\begin{eqnarray*}
\hat{\bu_1} & = & \bmat{c} \beta_1 (\hat{\mu}+\ab) \\  \beta_2 \hat{\mu} \\ \hat{\mu}(\hat{\mu}+\ab) \endbmat \\
& = & (1+\delta_1) \bu_1 + \bmat{c} 0 \\  0 \\ \mu^2(\delta_1+\delta_1^2) \endbmat
\end{eqnarray*}
and similarly
$$
\hat{\bu_3}  = (1+\delta_2)\bu_3 +  \bmat{c} 0 \\  0 \\ \nu^2(\delta_2+\delta_2^2) \endbmat
$$
This leads us to the useful fact that
\begin{eqnarray*}
\hat{\bu_1}^T \hat{\bu_3} & = & \nu(\nu+\ab)\mu^2(1+\delta_2)(\delta_1+\delta_1^2)
+ \mu(\mu+\ab)\nu^2(1+\delta_1)(\delta_2+\delta_2^2) 
+ \mu^2\nu^2(\delta_1+\delta_1^2)(\delta_2+\delta_2^2)\\
& = & \nu(\nu+\ab)\mu^2\delta_1 + \mu(\mu+\ab)\nu^2\delta_2 + O(\epsilon^2)
\end{eqnarray*}
Finally we note that $\Vert \bu_1\Vert \geq \mu(\mu+\ab)$ and $\Vert \bu_2\Vert \geq \nu(\nu+\ab)$ whence
$$
\frac{\hat{\bu_1}^T \hat{\bu_3}}{\Vert \bu_1\Vert \Vert \bu_2\Vert} \leq 
\frac{\mu^2}{\mu(\mu+\ab)}\delta_1 + \frac{\nu^2}{ \nu(\nu+\ab)}\delta_2 + O(\epsilon^2) \leq 2\epsilon + O(\epsilon^2)
$$
and the computed eigenvectors are numerically orthogonal provided that condition \ref{munu} is satisfied.

\subsection{Finding the Rightmost Eigenvalue of a Fully Reduced Arrow}
\label{algsec}

In this section we develop two stable and efficient methods for finding the rightmost eigenvalue of a fully reduced arrow matrix
$$
\Ab = 
\bmat{ccc}
0 & 0 & \beta_1 \\
0 & -\ab & \beta_2 \\
\beta_1 & \beta_2 & \gb \\
\endbmat
$$

Note that the block Gauss factorization of $\Ab -\l I$ is
$$
\Ab- \l I =
\bmat{ccc}
1 & 0 & 0 \\
0 & 1 & 0 \\
\frac{-\beta_1}{\l} & \frac{-\beta_2}{\l+\ab} & 1 \\
\endbmat
\bmat{ccc}
 - \lambda & 0 & \beta_1 \\
0 & -(\l + \ab) & \beta_2 \\
0 & 0 & -f(\l) \\
\endbmat
$$
where $f$, the {\em spectral function} of $\Ab$, is given by
$$
f(\lambda) = 
\l - \gb  - \frac{\beta_1^2}{\l} - \frac{\beta_2^2}{\l+\ab} .
$$
This is a rational Pick function with a pole at infinity. Inspection of the graph of the spectral function reveals that the elements of the shaft interlace the eigenvalues
\bne
\l_1 > 0 > \l_2 >  -\ab > \l_3 .
\label{interlace}
\ene
Moreover,  the derivative of the spectral function is
\bne
f^{\prime}(\l) = 1 + \frac{\beta_1^2}{\l^2} + \frac{\beta_2^2}{(\l + \ab)^2} .
\label{derspec}
\ene
and is clearly bounded below by one so that its zeros are, in a certain sense, well determined. Furthermore, we note that the second derivative of the spectral function
$$
f^{\prime\prime}(\l) = -2\left(\frac{\beta_1^2}{\l^3} + \frac{\beta_2^2}{(\l + \ab)^3}\right)
$$
is strictly negative over the interval $(0,\infty)$ and therefore $f^{\prime}(\l)$ is strictly decreasing over the same interval.

\subsubsection{Using the Borges-Gragg zero finder}

One approach to finding the unique zero of $f$ in the interval $(0 , +\infty)$ is to use the zero-finder developed in \cite{borgesgragg}.  Let $x_0$ be an initial approximation to the eigenvalue. If $x_j$ is known let our approximating function be 
$$
\phi_j (x) = \omega_0 x - \sigma - \frac{\omega_1}{x}.
$$
If we select the constants $\sigma$, $\omega_0$, and $\omega_1$ so that
\bne
\phi_j^{(i)}(x_j) = f^{(i)}(x_j) , \quad i = 0,1,2.
\label{rat_conds}
\ene
then we will be able to guarantee cubic convergence. Therefore, we solve
$$
\bmat{r c r}
-1 & x_j & -1/x_j \\
0 & 1 & 1/x_j^2\\
0 & 0 & -2/x_j^3
\endbmat
\bmat{c}
\sigma \\
\omega_0 \\
\omega_1
\endbmat
=
\bmat{c}
f(x_j) \\
f^{\prime}(x_j) \\
f^{\prime\prime}(x_j)
\endbmat
$$
and find
\begin{eqnarray*}
\omega_1 & = & \beta_1^2 + \beta_2^2 \left( \frac{x_j}{x_j + \ab} \right)^3, \\
\omega_0 & = & 1 +  \beta_2^2 \frac{ \ab}{(x_j + \ab)^3} , \\
\sigma & = & \omega_0 x_j - \frac{\omega_1}{x_j} - f(x_j).
\end{eqnarray*}
Note that $\omega_1 > \beta_1^2 > 0$ and also that  $\omega_0 > 1$. The inequalities are strict and since both $\omega_0 > 0$ and $\omega_1 > 0$ it follows that $\phi_j$ is a Pick function and has a unique zero $x_{j+1} \in (0 , +\infty)$. 

Casual inspection of the error function
$$
f(x) - \phi_j(x) = (1-\omega_0) x - (\gb - \sigma) -  \frac{\beta_1^2 - \omega_1}{x} - \frac{\beta_2^2}{x + \ab} 
$$
over the interval $(0,+\infty)$ reveals that it converges to $+\infty$ at the left boundary (as $x\rightarrow0^+$) and to $-\infty$ at the right boundary (as $x\rightarrow  +\infty$) and therefore crosses zero in the interval. Moreover, it is not hard to verify by differentiation and a bit of calculus that its derivative never changes sign in the interval. These facts imply that it crosses zero exactly once in the interval and further that we obtain monotonic convergence from any starting guess whatsoever $x_0 \in (0,+\infty)$. The cubic rate of convergence follows from (\ref{rat_conds}).

Successive iterates can be found by solving quadratic equations.  Rather than solve $\phi_j (x) = 0$ for $x_{j+1}$ it is better to solve
$$
\phi_j (x_j - \Delta) = 0
$$
for the {\em increment} $\Delta = x_j - x_{j+1}$.  Some rearrangement using (\ref{rat_conds}) reduces this to
\bne
a \Delta^2 + b \Delta - f = 0 ,
\label{del_quad}
\ene
with
\begin{eqnarray*}
a & = &  -\frac{\omega_0}{x_j} , \\
b & = &  f^{\prime} (x_j) + \frac{f(x_j)}{x_j} .
\end{eqnarray*}

 If $x_0 > \l_1$ then $f(\l) > 0 $ so $b > f^{\prime} > 1$ and $\Delta$ may be computed stably using
\bne
\Delta = \frac{2f/ b}{1 + \sqrt{1 + \frac{2a}{b}
\frac{2f}{b}}} ,
\label{increm_eq}
\ene
and we are therefore inclined to start to the right of $\l_1$. We can guarantee this by using the fact that the Borges-Gragg zero finder can start from $+\infty$ and choosing the first iterate from $+\infty$ to be our starting point. To do so, note that as $x \rightarrow +\infty$ the approximate Pick function tends to
\bne
\phi (x) = x - \gb - \frac{\beta_1^2+\beta_2^2}{x}.
\label{ext_start}
\ene
We propose to take $x_0$ to be the zero of (\ref{ext_start}) in $(0, +\infty )$  which is
$$
x_0 = \frac{\gb}{2} + \sqrt{ \left( \frac{\gb}{2} \right)^2 + \beta_1^2+\beta_2^2}.
$$
If $\gb<0$ then we may wish to multiply top and bottom by the conjugate to avoid cancellation.

Termination in this case is straightforward. The mean value theorem gives
$$
 \fp(c) = \frac{f(x_j) - f(\l_1)}{x_j - \l_1} = \frac{f(x_j)}{x_j - \l_1} 
$$
for some point $c\in(\l_1,x_j)$. However, since the iteration is monotonic our iterates are on the right so that $x_j > l_1$. Since $\fp(\l)$ is strictly decreasing over the same interval we can use $ \fp(x_j)$ as a lower bound on $ \fp(c)$
$$
 \fp(x_j)  < \frac{f(x_j) }{x_j - \l_1}
$$
and we conclude that the absolute error\footnote{The absolute value is unecessary here but we leave it in for pedagogical reasons.} is bounded by
$$
|x_j - \l_1| < \frac{f(x_j)}{ \fp(x_j) }.
$$
It is therefore reasonable to terminate when 
$$
\frac{f(x_j)}{ \fp(x_j) } < C \epsilon x_j
$$
for some constant $C$ since that will imply that $x_j$ has been found to high relative precision.

\subsubsection{Using Newton's Method}

It is remarkably easy to find the unique zero of $f$ in the interval $(0 , +\infty)$ using Newton's method. It is worth observing that because the derivative of the spectral function \ref{derspec} is strictly decreasing and bounded below by $1$ on the interval $(0 , +\infty)$ Newton's method will converge {\em monotonically} to the unique zero in the interval provided our initial guess, $x_0$, lies to the left of the zero. We can find an appropriate starting guess by finding the point where the spectral function intersects
$$
 -\frac{\beta_2^2}{x+\ab}
$$
since the graph of this function lies strictly below zero over the interval. After a bit of algebra and invoking the quadratic formula we find that this leads to a starting guess of
$$
x_0 = \frac{\gb}{2} + \sqrt{\left(\frac{\gb}{2}\right)^2 + \beta_1^2}.
$$

Termination in this case is also straigthforward. Since this iteration is monotonic on the left we use the fact that $f^{\prime}(\l) > 1$ over the interval $(x_j,\l_1)$ and we conclude that the absolute error is bounded by $ |x_j - \l_1| < | f(x_j)|.$ It is therefore reasonable to terminate when $ |f(x_j)| < C \epsilon x_j$ for some constant $C$ since that will imply that $x_j$ has been found to high relative precision.

\subsection{Testing}

In this section we test the algorithm we have developed herein, which we call {\tt threig}, versus the \matlab function {\tt eig} which uses routines from \textsc{LAPACK}. Our test is designed to compare the performance of our algorithm to that of the industry standard and therefore the protocol is simple. 
\begin{enumerate}
\item
Create a set of 100,000 random real symmetric test matrices where the elements of each test matrix $T$ are {\em randomly} drawn from a known distribution (the distributions we test are uniform, normal, and chi-square).
\item
Run both {\tt threig} and {\tt eig} on each test matrix $T$.
\item
Compute two measures of accuracy for each case:
\begin{enumerate}
\item
$\Vert I - V^TV \Vert_F$ which measures the orthogonality of the computed eigenvectors.
\item
$\Vert TV - V\Lambda \Vert_F$ which measures the accuracy of the computed spectral factorization (i.e. how well it reconstructs the test matrix.)
\end{enumerate}
\item
Finally we compute difference between the results from {\tt eig} and those from {\tt threig}, so that the ouptut will be negative when the measured error in {\tt threig} exceeds that of {\tt eig} and positive otherwise. We then sort these values and generate a plot of the results. Since the results are sorted before being displayed the plot reveals how often and by how much one algorithm is {\em superior} to the other.
\end{enumerate}

\begin{figure}[ ht ]
\begin{center}
\includegraphics[width=6in]{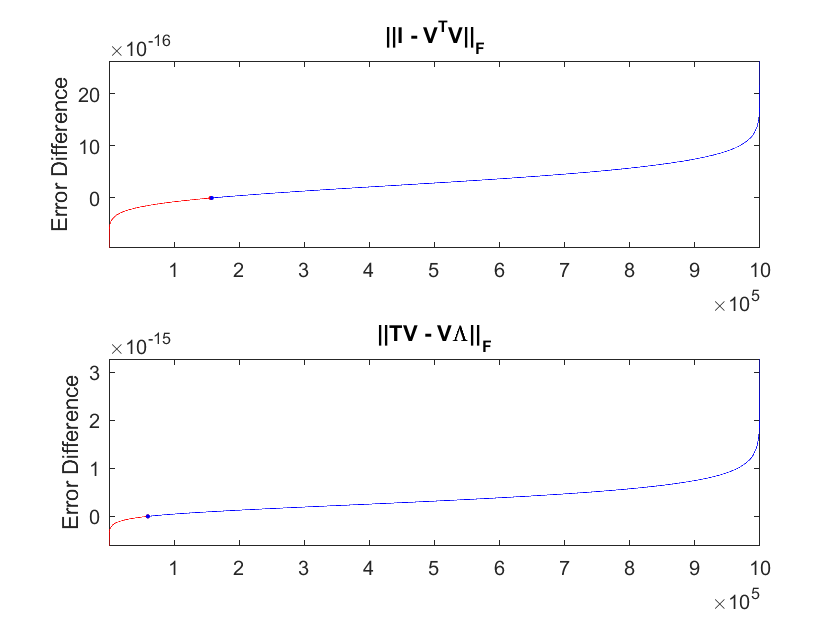}
\caption{Test matrix elements distributed $U(0,1)$.}
\end{center}
\end{figure}

\begin{figure}[ ht ]
\begin{center}
\includegraphics[width=6in]{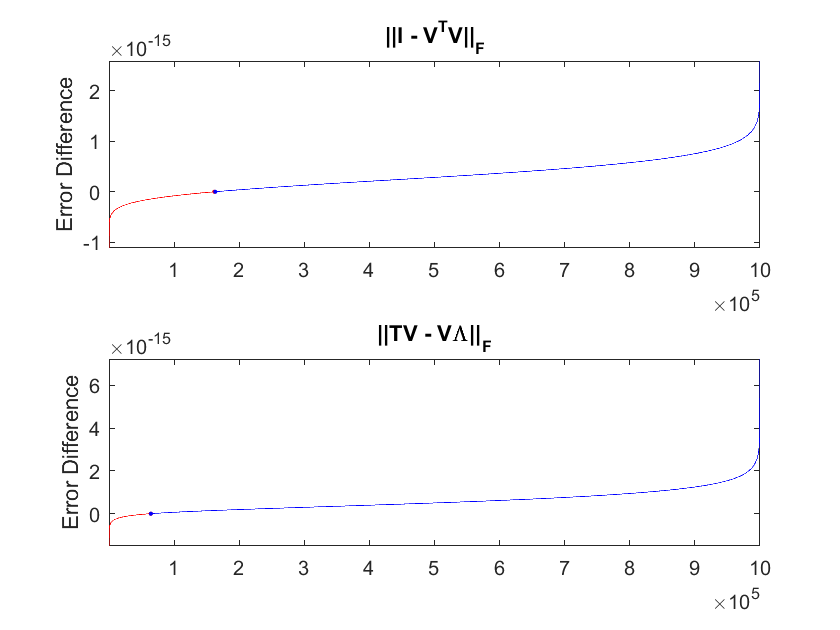}
\caption{Test matrix elements distributed ${\cal N}(0,1)$.}
\end{center}
\end{figure}

\begin{figure}[ ht ]
\begin{center}
\includegraphics[width=6in]{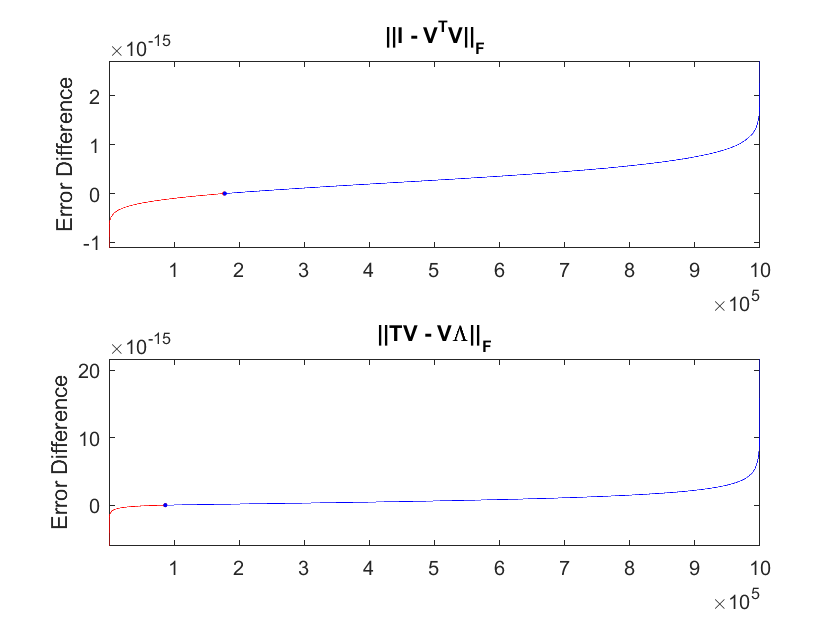}
\caption{Test matrix elements distributed $\chi^2_1$.}
\end{center}
\end{figure}

On viewing the figures we see that the proposed algorithm is consistently more accurate in both measures.

\end{document}